 
%
\input amstex

\def\aaaa#1#2#3{A^{#1}_{{#2}{#3}}}

\def\bcc#1{\Bbb C^{#1}}
\def\bpp#1{\Bbb P^{#1}}
\def\brr#1{\Bbb R^{#1}}
\def\bee#1{\Bbb E^{#1}}

\def\cf{\Cal F}
\def\ci{\Cal I}

\def\ee#1{e_{#1}}

\def\frp#1#2{\frac{\partial {#1}}{\partial {#2}}}

\def\hhh#1#2#3{h^{#1}_{{#2}{#3}}}
\def\hd{, \hdots ,}

\def\inv{{}^{-1}}

\def\ns{n+s}
\def\nr{n+r}
\def\na{n+a}
\def\ooo#1#2{\omega^{#1}_{#2}}
\def\oo#1{\omega^{#1}}

\def\ot{\!\otimes\!}
 
 \def\pii#1{\pi^{#1}}
\def\upp#1#2{p^{#1}_{#2}}

\def\rr#1#2#3#4{R^{#1}_{{#2} {#3}{#4}}}
 \def\ra{\rightarrow}

\def\rdiff{\}_{diff}}

\def\ta#1{\theta^{#1}}
\def\taa#1#2{\theta^{#1}_{#2}}
\def\tor#1#2#3{T^{#1}_{{#2}{#3}}}

\def\ttrace{\text{trace}}
\def\tdim{\text{dim}}

\def\tspan{\text{span}}
\def\tmod{\text{ mod }}

\def\pii#1#2{\pi^{#1}_{#2}}

\def\uu#1#2{ u^{#1}_{#2}}

\def\up#1{{}^{({#1})}}
\def\upperp{{}^{\perp}}

\def\ww{\wedge}
\def\xx#1{x^{#1}}

\def\na{n+a}

\define\intprod{\mathbin{\hbox{\vrule height .5pt width 3.5pt depth 0pt %
        \vrule height 6pt width .5pt depth 0pt}}}
\documentstyle{amsppt}
\magnification = 1000
\hsize =14truecm
\hcorrection{.5truein}
\baselineskip =18truept
\vsize =22truecm
\NoBlackBoxes
\topmatter
\title Exterior differential systems: a geometric approach to PDE\linebreak
Lecture notes from the 1997 Daewoo workshop
\endtitle
\rightheadtext{Exterior Differential Systems}
\author
J.M. Landsberg
\endauthor

\date {  July 7-11, 1997}\enddate

\address{ Laboratoire de Math\'ematiques,
  Universit\'e Paul Sabatier, UFR-MIG
  31062 Toulouse Cedex 4,
  FRANCE}\endaddress
\email {jml\@picard.ups-tlse.fr }
\endemail

\endtopmatter

\document

\heading Lecture 1: Introduction to exterior differential systems.
\endheading

In this lecture we will see how to study a system of
partial differential equations (henceforth denoted pde) from a geometric
perspective.

\smallpagebreak

\noindent{\bf Example 1, the minimal surface equation}. Given a surface 
in Euclidean $3$-space, $M^2\subset\bee 3$, 
described locally as a graph $z=z(x,y)$, at a point
$p=(x,y,z(x,y))$, let
$$
H(p)=   \frac{1}{2}\frac{(1+z_y^2)z_{xx} - 2z_xz_yz_{xy} + (1+z_x^2)z_{yy}}
{(1+z_x^2 + z_y^2)^\frac{3}{2}}.  
$$
$H(p)$ is called the {\it mean curvature function} of the surface $M$.
$H(p)$ has geometric meaning, that is, it is a well defined function
on $M$, independent of coordinates chosen and invariant  under
the displacement of $M$ by
rigid motions.

Consider the pde for surfaces in $\bee 3$ with the property that
$H(p)\equiv 0$. Geometrically such surfaces are critical for the
variation of area and are called {\it minimal surfaces}. If $M$ is a
minimal surface, then
for all $x\in M$, there exists a open neighborhood $\tilde U_x\subset M$
such that for all open subsets
$U\subset \tilde U_x$ and patches of surface $V$ with
$\partial V=\partial U$, then $\text{area}(V)\geq\text{area}(U)$
(see e.g.  [S, III]). 
$$
\
$$
\bigpagebreak

$$
\
$$
\bigpagebreak

$$
\
$$

 We will be studying minimal surfaces
and higher dimensional generalizations in the next lecture.
One amazing aspect of the pde for minimal surfaces
in $\bee 3$ is that the solutions are exactly given
by solutions to the Cauchy-Riemann equations. More precisely,
there is  the following  
theorem:

\proclaim{Theorem (Weierstrass [S, IV, p 395])} Minimal surfaces
$M^2\subset\bee 3$ can be described by holomorphic functions.
More precisely, every minimal surface in $\bee 3$
 is locally of the form:
$$
\align
&x(w,\overline w)= Re\frac 12\int f(w)(1-g(w)^2)dw\\ 
&y(w,\overline w)= Re\frac i2\int f(w)(1+g(w)^2)dw\\ 
&z(w,\overline w)= Re\int f(w)g(w)dw\endalign
$$
where  $g(w)$ is a meromorphic function
and   $f(w)$ is a holomorphic function vanishing precisely at
the poles of $g$ and the integrals are path integrals beginning at
some fixed point $w_0$ to the point $w$.
\endproclaim

There is actually ambiguity in the formulae and, as we will see, the
minimal surfaces essentially depend on one holomorphic function.

How could
one discover such a beautiful theorem? 
We will see that by rephrasing the minimal surface equation
as an exterior differential system, one might guess that such
a formula exists, and in fact be able to see the geometry of such.
(Of course Weierstrass had to work alot harder since he didn't have
the machinery available to him.) The exterior differential systems
machinery has been used to recognize other systems of pde
appearing in geometry as the Cauchy-Riemann system and other
familiar systems.

The minimal surfaces $M^2
\subset\bee{2+s}$, are also all described in terms of holomorphic functions
but this fails to be true   when $\tdim M>2$.

On the other hand, one might hope that special classes of minimal
submanifolds might still be described by holomorphic data.
(This is because the minimal submanifold equation \lq\lq looks
like\rq\rq\ Laplace's equation (they have the same symbol) so just
as in any even dimension one can find special solutions of Laplace's
equation by solving the Cauchy-Riemann equations (or more generally
solutions of a Dirac equation) one could hope to do the same for
minimal submanifolds. When we have understood goal 1 below, we
will see how to arrive at the description of minimal surfaces
in terms of holomorphic data and how to look for special classes
of minimal submanifolds of arbitrary dimension.

\medpagebreak

We will have the following general goals for these lectures:

\smallpagebreak

\noindent {\it 1. To explain how to find an appropriate geometric 
setting for studying a given system of pde.}

 To do this we will first work on a larger space to eliminate choices
of coordinates that have no geometric meaning, and then 
we will \lq\lq quotient out\rq\rq\ 
to a smaller space where variables not relevant to the system are
eliminated.

 On the geometrically determined space one can then hope to recognize
familiar systems in disguise, and to find special classes of solutions
(e.g. solutions with \lq\lq symmetry\rq\rq ).

\smallpagebreak

\noindent{\it 2. To explain the  Cartan algorithm to determine the moduli
space of local solutions to any given exterior differential
system and an appropriate initial value
problem for the system. }

The general Cartan algorithm is somewhat complicated to state, so for
the moment, I'll just
mention a theorem we will study in lecture 3 that is proved using the
algorithm:

\proclaim{ BCJS  Theorem  (see [BCG${}^3$, p 302],[S, V, p 216])}
Let $g$ be any analytic Riemannian metric in a neighborhood $V$ 
of $0\in\brr n$. 
There exist local analytic isometric embeddings
on some smaller neighborhood $U\subseteq V$,
$i: (U,g)\ra \bee{ \binom{n+1}2}$. The embeddings depend on $n $ functions
of $n-1 $ variables.
\endproclaim

Note the form of the conclusion of the BCJS  theorem. This is the type of
answer the Cartan-K\"ahler theorem provides, not only an assertion
of existence, but a rough description of the size of the moduli space
of solutions.

Recently I have utilized the Cartan machinery to study questions
in algebraic geometry. It has been quite useful  for sorting
out what of the geometry of projective varieties arises from
global considerations and what is due to the local projective
geometry. (See [L, 3-8]).

\bigpagebreak

\noindent{\bf Linear Algebra Aside}. Let $V$ be a vector space ($\brr n$ or
$\bcc n$
for our purposes), and let $V^*$ denote the dual vector space.
Let $\ta 1\hd\ta s\in V^*$ be linearly independent.
One has the correspondence
$$
\matrix
s-\text{planes in }V^* & \leftrightarrow & (n-s)-\text{planes in }V\\
\tspan\{\ta 1\hd \ta s\}& &
   \{ v\in V | \ta a(v)=0, 1\leq a\leq s\}\endmatrix
$$

\smallpagebreak

\noindent{\bf Example 2}. Consider the system of pde  
$$
u_x = A(x,y,u) \ \ \ u_y = B(x,y,u) \tag 2.1
$$
for one function of two variables
(where $A,B$ are given functions).
Note that it is overdetermined (two equations for one function) so
one expects that there won't be any solutions unless $A$ and $B$
satisfy a compatibility condition. To phrase this system geometrically, let
$\brr 3$ have coordinates $(x,y,u)$ and consider the forms
$$
\ta{} := du -Adx-Bdy \ \ \ \Omega := dx\ww dy   
$$
Observe that
we have the correspondence

$$
\matrix 
 \text {\it Solutions to (2.1)}&\leftrightarrow&
\text{\it surfaces }i: S\ra \brr 3,
\text{\it such that   }\\
& & i^*(\theta)\equiv 0\text{ \it and }i^*(\Omega )\text{ \it is
nonvanishing}.
\endmatrix
$$
When using the exterior differential systems
machinery, the problem of finding
solutions to a system of pde is always  rephrased  
as a problem in submanifold geometry as in this example.

Geometrically, at each $p\in \brr 3$
there exists a unique two-plane $U_p\subset T_p\brr 3$ 
annhilated by $\theta$ so
if there exists a solution passing through $p$, its tangent
plane must be $U_p$. So the question is:
Do these two-planes \lq\lq fit together\rq\rq ?
 This is a second order condition. To determine the condition
in coordinates, we utilize   that
  mixed partials commute:
$$
\align 
(u_x)_y  & = \frp {} y A(x,y,u) = A_y(x,y,u) + \frp u y A_u(x,y,u) = A_y
+BA_u  \\
(u_y)_x & = B_x +AB_u \endalign
$$
Setting $(u_x)_y = (u_y)_x$  uncovers the equation
$$
A_y + BA_u = B_x + AB_u \tag 2.2
$$

More geometrically, $(u_x)_y= (u_y)_x$ if and only if
$$
d\theta \equiv 0\tmod \{\theta\}\tag 2.3
$$

\noindent{\bf Notation}:
If $v_1\hd v_s$ are vectors, $\{ v_1\hd v_s \}$
denotes their linear span. If
$\ta 1\hd \ta s\in\Omega^*(M)$ are forms
of any degree, let  $\ci=\{\ta 1\hd\ta s\}_{alg}$
denote the algbraic ideal generated by $\ta 1\hd \ta s$. Given
$\beta\in \Omega^k (M)$, we say $ \beta\equiv 0\tmod \{\ta 1\hd\ta
s\}_{alg}$
 if
$\beta = \alpha_1\ww\ta 1 +\hdots + \alpha_s\ww\ta s$ for some
$\alpha_1\hd\alpha_s$.

\smallpagebreak

Clearly (2.3) is a necessary condition, but in fact it is sufficient.

\proclaim{Frobenius Theorem (see [BCG${}^3$, p 27])}   Let $\ta 1\hd\ta s$
be a pointwise linearly
independent one-forms on a manifold $\Sigma^{n+s}$ of dimension $n+s$.
If for all $x\in \Sigma$ and for all $1\leq a\leq s$, we have
$$
d\ta a\equiv 0\tmod \{\ta 1\hd \ta s \} \tag 2.4
$$
then through each $x\in \Sigma$ there exists a unique submanifold
$i: M\ra \Sigma$ such that for all $y\in M$, $T_yM= \{\taa 1y\hd\taa
sy\}\upperp$.
In this case we say the distribution $ \{\ta 1\hd \ta s\}\upperp$ is
{\it integrable}.\linebreak
 ($\upperp$ denotes the annhilator in the dual vector space.)
\endproclaim

For example 2, there are the following cases:

Case 1. The system is integrable.

Case 2. Otherwise, to find a solution, restrict to $\Sigma '\subset\brr 3$,
 the surface
where (2.2) holds. In the neighborhood of a $p\in \Sigma '$ there are two 
possibilities. 

Case 2a. Since $\tdim\Sigma ' =2$, we already have a surface so
if $u$ appears in  the equation defining $\Sigma '$, 
equivalently, if $i^*(\Omega )$ is nonvanishing, then solving for $u$
gives  a candidate
for a solution and we need to check if this function satisfies the system,
equivalently, if $i^*(\theta)=0$.

Case 2b. If $u$ does not appear, then on $\Sigma '$ there is a 
relation between our independent variables $x$ and $y$ so
$i^*(\Omega )=0$
and there are no solutions.  
\smallpagebreak

\noindent{\bf Definitions}:
 An {\it exterior differential system with independence
condition} on a manifold $\Sigma$ consists of an ideal
 $\ci\subset\Omega^*(\Sigma)$ closed
 under exterior differentation and a differential $n$-form
$\Omega\in\Omega^n(\Sigma)$. ($\Omega^*(\Sigma) = 
\oplus \Omega^k(\Sigma)$ denotes the space of all differential forms on
$\Sigma$.)

An {\it integral manifold}
(solution) of the system $(\ci , \Omega )$ is an immersed $n$-fold
$f:M^n\ra \Sigma$ such that $f^*(\alpha ) = 0\ \forall \alpha\in\ci$ and
$f^*(\Omega )$ is nonvanishing.

For example, in example 2, $\ci =\{ \theta\}_{diff}= \{ \theta ,d\theta
\}_{alg}$.
One could rephrase (2.4) as saying if 
$\{\ta 1\hd \ta s\}_{alg}=\{\ta 1\hd\ta s\}_{diff}$, then the system
is integrable.

It is useful to
study integral manifolds infinitesimally.
Let  $V$ be a vector space, and let $G(n,V)$ denote the 
Grassmanian of $n$-planes through the origin in 
$V$.  We define an {\it integral element} of $(\ci ,\Omega)$
at $x\in\Sigma$
 to be an $E\in G(n,T_x\Sigma )$ such that
 $ \alpha | _E = 0\ \forall \alpha\in\ci$ and
$ \Omega |_E\neq 0$.
Integral elements are the candidates for tangent spaces to integral
manifolds.

In these lectures we will actually only study a special class of exterior
differential
systems with independence condition. A {\it Pfaffian} eds is one generated
by differentially one-forms, say $\ta a$, $1\leq a\leq s$.
I.e. $\ci = \{ \ta a\}_{diff}=\{ \ta a, d\ta a \}_{alg}$.
 Let  $I_x:=\{\ta a|_x\}
\subset T^*_x\Sigma $
be the distribution they generate. Write the indpendence condition
as $\Omega = \oo 1\ww\hdots\ww \oo n$, and let $J_x:=\{\ta a|_x, \oo
i|_x\}$.
The Pfaffian system is {\it linear} if 
$$
d\ta a\equiv 0\tmod J.
$$
We will restrict our attention to linear Pfaffian systems.
The name \lq\lq linear\rq\rq\ is used because if a Pfaffian system
is linear, then its integral elements at a point are determined by linear
equations, see e.g. example 3 below.
One often denotes a linear Pfaffian system by $(I,J)$ instead
of $(\ci,\Omega)$. 
\smallpagebreak

 In principle we can study
any exterior differential system with independence condition 
just by studying linear Pfaffian systems
as follows:

\subheading{The canonical linear Pfaffian system on $G(n,T\Sigma)$}
Let $\Sigma$ be any manifold.
Consider $\pi : G(n,T\Sigma)\ra \Sigma$,
 the Grassmann bundle of all $n$-planes in all tangent spaces to $\Sigma$. 
 We will denote points of $G(n,T\Sigma)$ by $(p,E)$ where $p\in \Sigma $
 and $E\subset T_p\Sigma$ is an $n$-dimensional subspace.  $G(n,T\Sigma)$
 carries a cannonical linear Pfaffian system $(I,J)$  on it.
  Namely it is the  system whose solutions are exactly the lifts to
 $G(n,T\Sigma)$  of mappings $f: X^n \ra M$, where the lift of $f$ is 
$\tilde f(x) = (f(x), T_{f(x)}f(X))$.  The system is defined by
$$
I_{(p,E)} := \pi^*(E^{\perp}) \ \ \ J_{(p,E)} := \pi^*(T^*_p\Sigma)
$$

\subheading{The prolongation of an exterior differential system}
Now say $\Sigma$ is equipped with an exterior differential system
with independence condition $(\ci,\Omega)$.
Let 
$$
\tilde\Sigma:=
\{ (x,E) \in G(n,T\Sigma ) \mid E\in G(n, T_x\Sigma )
\text{ is an integral element of }
(\ci,\Omega)\text{ at }x\}. 
$$

The pullback  of the canonical system $(I,J)$ on
$G(n,T\Sigma )$ to $\tilde \Sigma$
 is 
a linear Pfaffian system called the {\it prolongation} of $(\ci,\Omega)$.

\smallpagebreak

 Any first order system of pde for
maps $u: \brr n\ra \brr s$ may be described as a
linear Pfaffian system as follows:

\subheading{The canonical contact system on the space of one-jets}
Let $J^1=J^1(\brr n ,\brr s ) := {\brr n} \times {\brr s} \times
{\brr {ns}}$ with coordinates
$(\xx i,u^a,p^a_i)$. $J^1$ is called the {\it
space of one-jets of mappings from $\brr n$ to $\brr s$}.
 On
$J^1$ let 
$$
\align &
\ta a := du^a - p^a_i dx^i, \ \  1\leq a\leq s\tag  j1.1\\
&\Omega := dx^1\ww\hdots\ww
dx^n.\endalign
$$
Integral manifolds of the system $(I:=\{\ta a\} , J:=\{\ta a,  dx^i\})$
 project to graphs of maps $u:\brr n\ra \brr s$. 
We call the system $(I  , J )$
the {\it canonical contact system on $J^1$}.
 One can also work in
the complex category with holomorphic functions and then we will
call such a system a {\it complex contact system}.

\smallpagebreak

Note that integral manifolds of (j1.1)
  are describable in terms of $s$ arbitrary
functions of $n$ variables. Namely
choose functions $f^a(\xx 1\hd\xx n )$, $1\leq a\leq s$ and set
$$
\align u^a&= f^a(\xx 1\hd\xx n)\\
\upp ai &= \frp {f^a}{\xx i}(\xx 1\hd\xx n). 
\endalign
$$
  If one is ever
so lucky to recognize a given system of pde as a contact system in
disguise, one can simply write down local solutions in terms of arbitrary 
functions.

\smallpagebreak

\subheading{How to express a system of pde as a linear Pfaffian system}
Now say we have a first order system of pde:
$$
F^r (x^i,u^a, \frp {u^a}{x^i} ) =0,  \ \ 1\leq r\leq R  \tag  j1.2
$$
($R$ equations). Consider the codimension $R$ submanifold 
 $\Sigma\subset J^1$ given
by $F^r(x^i,u^a, p^a_i ) =0,\ 1\leq r\leq R$.
  Solutions to (j1.2) correspond to
integral manifolds of the restriction of the
contact system on $J^1$ to $\Sigma$.  

\smallpagebreak

Note that lifting   up to $J^1$ has the effect of turning
the derivatives into independent variables. Restricting to $\Sigma$
has the effect  of forcing the pde to hold as an algebraic (or analytic)
relation among the variables. The    integral submanifolds
of $\Sigma$ are those where the variables representing the
derivatives behave like derivatives and the independent variables
stay independent.
 
\smallpagebreak

\smallpagebreak

\noindent{\bf Example 3, the Cauchy-Riemann equations}. 
Write the Cauchy-Riemann equations as
$\uu 1{\xx 1}=\uu 2{\xx 2},\uu 1{\xx 2}=-\uu 2{\xx 1}$.
Let  $J^1(\brr 2,\brr 2)$ have 
coordinates $(x^i, u^a, p^a_i)$, $1\leq i,j\leq 2, 1\leq a,b\leq 2$,
let $\Sigma^6\subset J^1$ be defined
 by the equations $\upp 11-\upp 22=0, \upp 12 + \upp 21 =0$.
The  restriction of the
canonical contact system 
$(I,J)$ on $J^1$ to  $\Sigma$ is equivalent to the Cauchy-Riemann
equations.

Following  
  example 2, and having no idea what else to do, we differentiate.
We have
$$\align
d\pmatrix \ta  1\\ \ta 2\endpmatrix  &= 
\pmatrix -d\upp 11 & - d\upp 12 \\ - d\upp 21 & - d\upp 22\endpmatrix \ww
\pmatrix d\xx 1 \\ d\xx 2\endpmatrix  \\
& = \pmatrix  \pi_1 & \pi_2 \\  -\pi_2 & \pi_1\endpmatrix
\ww
\pmatrix d\xx 1 \\ d\xx 2\endpmatrix  .\endalign
$$
where we use the notation $\pi_1 =-d\upp 11$, $\pi_2 = -d\upp 12$.

Note that $\ta 1,\ta 2,\pi_1,\pi_2,d\xx 1,d\xx 2$ form
a basis of $T^*_x\Sigma$ for all $x\in\Sigma$ and the space of integral
elements at   $x$ is given by the  linear equations:
$$
\align
\ta 1&=0\\
\ta 2&=0\\
\pi_1 &= ad\xx 1 + bd\xx 2\\
\pi_2 &= bd\xx 1 -ad\xx 2\endalign
$$
where $a,b$ are constants. Thus there is a two dimensional space
of integral elements through each point, in contrast to the case of
the Frobenius theorem where there was a unique integral element passing
through each point.

Since $\Sigma$ is a linear space of even dimension it can be given a
complex structure,
but what is important here is that it can be given a complex structure such
that $\ci$ becomes
a  {\it complex contact structure}.

On $T^{*\Bbb C}\Sigma$, let
$$
\align
\theta & := \ta 1 + i\ta 2   \\
\omega &:= d\xx 1 + id\xx 2 \\
\pi &:= \pi_1 + i\pi_2\endalign
$$
define an
almost complex structure
on $\Sigma$ by specifiying
a basis of $T^{*10}\Sigma$. This almost complex structure is integrable.
Moreover the system can be written as
$\ci = \{ \theta , \overline\theta \rdiff$, $\Omega =
\omega\ww\overline\omega$.

$(\ci,\Omega )$ is a  {\it complex contact system} in the sense that there exist
complex coordinates
$(w, p, z)$  on $\Sigma$ such that 
$\ci = \{\theta := dw - pdz,\overline\theta :=
d\overline w - \overline pd\overline z\}$ and integral manifolds are
obtained by letting $f(z)$ be any holomorphic function
and setting 
$$
w = f(z),\  p = f'(z).
$$
 
\medpagebreak

 In the next lecture we will see how to recognize the eds for minimal
surfaces as being equivalent to the Cauchy-Riemann equations.
We will assume some familiarity with moving frames so those of you
not yet familiar with moving frames are encouraged to 
do the preparatory exercises at the end of this lecture.

So far we have seen how to enlarge the space we work on, as stated
in goal 1. Now I would like to describe how to \lq\lq 
quotient out\rq\rq\
to a smaller space.

\smallpagebreak

\noindent{\bf Definition}. Let $\ci$ be and exterior differential
system on a manifold $\Sigma$. A vector field $\xi\in\Gamma (TM)$ is
said to be a {\it Cauchy characteristic vector field} if
$\xi\intprod \alpha =0$ for all $\alpha \in \ci$.

\proclaim{Proposition [BCG${}^3$, p 21]} Let $\xi_1\hd\xi_k$
be Cauchy characteristic vector fields, then the distribution
they span is integrable.
\endproclaim

Thus if $\ci$  on $\Sigma$
has Cauchy characteristic vector fields 
then we can actually work on a smaller space,
\lq\lq$\Sigma/\text{leaves}$\rq\rq . We will see several examples
of systems with Cauchy characteristic vector fields in the next
two lectures.

\medpagebreak

In the few minutes remaining, I would like to mention 
a few results I obtained recently using the Cartan algorithm.
Motivated by some problems in algebraic geometry, I addressed
the following question:

{\it How many derivatives does one need to take to see if
a submanifold of affine or projective space is built out of linear
spaces?}

For example, given a surface, even in Euclidean $3$-space, if it is
negatively curved, then  at each point there will be two lines osculating
to order two (called the {\it asymptotic lines}). Thus to see if a surface
in $3$-space is ruled by lines one needs at least three derivatives.
It is a classical result that in fact three derivatives are enough.
Here is a generalization of this result:

\proclaim{Theorem [L7, 1]}
Let $X^n\subset\Bbb A^{\na}$ or $X^n\subset\bpp\na$ be a patch
of a smooth submanifold of an affine or projective space
such that at every point there is a line osculating to order $n+1$.
Then $X$ is ruled by lines.

  There exist patches $X^n\subset\Bbb A^{n+1}$ or $X^n\subset\bpp{n+1}$ 
having a line osculating to order $n$ at every point that are not ruled.
In fact over $\Bbb C$, every hypersurface has this property.   
\endproclaim

Given integers $(n,a,k)$ one could ask what is the smallest number
$m_0$ such that if $X^n\subset\bpp\na$ is a patch such that at each
point there is a $k$ plane osculating to order $m_0$, then
the $k$-planes are (locally) contained in $X$, yet there exist
examples of patches $X^n\subset\bpp\na$ with $k$-planes osculating
to orders $m_0-1$ that are not contained. Here is a partial result
regarding this question:

\proclaim{ Expectation Theorem, [L7, 6]} Let $(n,a,k,m)$ be   natural
numbers satisfying $m\geq 3$ and
$$
a[\binom{k+m-1}{m-1} -k -1]\geq k(n-k).
$$
Let $X^n\subset\Bbb A^{\na}$ or $X^n\subset\bpp\na$ be a patch
of a smooth submanifold of an affine or projective space
having the properties that at each $x\in X$ there exists a
$k$-dimensional linear space $L_x$,
disjoint from the fiber of the Gauss map,
 osculating to order $m$
and such a certain rank condition holds through level $m-1$.
 Then these linear spaces are locally contained in $X$.
  If $F$ denotes
the fiber of the Gauss map at $x$ and $\tdim L_x\cap F=\lambda$, then the
same
conclusion holds as above with $k,n$ respectively replaced by
$k-\lambda, n-\lambda$.
In fact, if the Gauss map is degenerate, one
can replace $L$ by the span of $L$ and $F$.
\endproclaim

The rank condition essentially requires that the differential invariants
of $X$ at $x$ are reasonably generic
among all possible tensors of differential
invariants allowing a $k$-plane to osculate to order $m$,
see [L7] for details.

I do not know if the $m$ furnished by the expectation theorem
is optimal. I showed it gave a lower bound when $k=1$ and was optimal
when $k=1$ and $a=1$ and that it gave an upper bound when $k=n-2$.
I am currently working on this question.

\subheading{Problems for lecture 1}
 
1a. Calculate the dimension of the space of integral elements
at a point of $J^1$
for the canonical contact system on $J^1$.

1b. Calculate the dimension of the space of integral elements
at a point of $\Sigma\subset J^1(\brr 2,\brr 1)$ where $\Sigma$
corresponds to the hypersurface induced by the equation $u_x+u_y=0$.
 
1c.  What can one say about
the dimension of the space of integral elements at a point
of a hypersurface $\Sigma\subset J^1$ (where one
restricts the canonical system to the hypersurface)?

2. Describe the  
construction
corresponding to $J^1$
 for second order systems of pde, then systems of $k$-th order.
In each case, calculate the dimension of the space of integral
elements of the canonical contact system.

\subheading{Problems to prepare for lectures 2 and 3}

1. Let $G\subset Gl(n,\brr{})$ be a matrix Lie group. Consider
  $g\in G$ as a map $g: G\ra n\times n\text{ matrices}$. Let
$g\inv$ be the inverse matrix of $g$. Let $\Omega:= g\inv dg$
be the corresponding matrix valued one-form. $\Omega$ is called
the {\it Maurer-Cartan} form of $G$. Show that in fact it is
$\frak g$ valued, where $\frak g$ is the Lie algebra of $G$.
Show $\Omega$ satisfies the {\it Maurer-Cartan equation},
$d\Omega = -\Omega\ww\Omega$.
(Hint: calculate $0=d(Id)=d(g\inv g)$.)

2. Prove the  {\it Cartan lemma}:
Let
$v_1\hd v_k$  be linearly independent elements of
a vector space $V$ and 
let $w_1\hd w_k$ be elements of $V$ such that
$w_1\ww v_1 + \hdots w_k\ww v_k=0$. Then there exist constants
$h_{ij}=h_{ji}, 1\leq i,j\leq k$ such that 
$w_i =\Sigma_jh_{ij}v_j$.

3. Let $(M^n,g)$ be a Riemannian manifold. Let $\cf\ra M$ denote
the bundle of all orthonormal frames. Write $f\in\cf$ as
$f=(x,\ee 1\hd\ee n)$,
where $x\in M$ and $\ee 1\hd \ee n$ 
is an orthonormal basis of $T_xM$. Show that we may write $dx=\oo i\ee i$,
$d\oo i= \ooo ij\ee j$ 
where $\oo i,\ooo ij\in\Omega^1 (\cf )$,
with $\ooo ij+\ooo ji=0$. The forms
$\ooo ij$ are called connection forms. In fact
$\nabla_{\ee i}X= dX(\ee i) + \ooo ji(X)\ee j$, where $\nabla$ is
the Levi-Civita connection. Show that we have the equations
$$
\align
d\oo i &=-\ooo ij\ww\oo j\\
d\ooo ij &= -\ooo ik\ww\ooo kj + \rr ijkl\oo l\endalign
$$
where
the $\rr ijkl$ are functions defined on $\cf$ such that
 $\rr ijkl\ee i\ot\oo j\ot \oo k\ot \oo l$ descends to
be a well defined element of $\Gamma (TM\ot T^{*\ot 3}M)$
which is the Riemann curvature tensor.

4. Let $V$ be a vector space, let $G(n,V)$ denote the Grassmanian
of $n$-planes through the origin in $V$. Show that for all $E\in G(n,V)$,
the tangent space to $E$ has the additional structure as a vector
space of linear maps, more precisely, that $T_EG(n,V)\simeq E^*\ot V/E$.
Hint: Write $E=v_1\ww\hdots\ww v_n$, consider a curve
$E(t)=v_1(t)\ww\hdots\ww v_n(t)$ and differentiate at $t=0$.
If $V$ comes equipped with an inner product, we may identify
$V/E$ with $E\upperp\subset V$.
 
\heading Lecture 2: Applications to the study of minimal submanifolds.
\endheading
 
\subheading{The frame bundle}
It is not usually possible to
choose coordinates adapted everywhere to the geometry
of a given problem, but it is usually
possible to choose adapted frames.

Let $V=\brr \ns$ have the standard
inner product (we write $V$ to emphasize the structure of a vector
space to avoid confusion with
Euclidean  space $\bee\ns$.)
Let $\cf_{\bee \ns}$ be the bundle of all 
orthonormal framings of $\bee \ns$,
that is each $f\in \cf$ is $f=(x,u)$, where $x\in \bee \ns$
and $u: T_p\bee \ns \ra V$. is an isometry.
 
Given a submanifold $M^n\subset\bee \ns$, we can 
define a subbundle of the restriction of $\cf $ to
$M$, namely, writing
$u= (\ee 1\hd\ee n,\ee{n+1}\hd\ee\ns )$, we may restrict further to frames such
that $T_xM=\{ \ee 1\hd \ee n\}$. We   denote this bundle
by $\cf^1\ra M$.  

The main advantage of frame bundles is that they come equipped
with a canonical framing, that is a canonical basis of their cotangent
space at each point.

Use index ranges $1\leq i,j,k\leq n,\ n+1\leq a,b\leq s$.
On $\cf\ra\bee\ns$ we may write:
$$
d(x,\ee 1\hd\ee\ns ) =(x,\ee 1\hd\ee\ns ) 
\pmatrix 0&0&0\\ \oo i &\ooo ij & \ooo ib\\
\oo a&\ooo aj &\ooo ab\endpmatrix
$$
where $\oo\alpha, \ooo\alpha\beta=-\ooo\beta\alpha\in \Omega^1(\cf)$.

On $\cf^1\ra M$, we have $\oo a=0$, which of course implies
that $d\oo a=0$.
(Here and in what follows we commit a standard abuse of notation, we
really mean the pullback of $\oo a$ to $\cf^1$ is zero.)
Since $\cf $ is a Lie group,  letting
$\Omega$ denote its Maurer-Cartan form, we have 
the Maurer-Cartan equation, $d\Omega =-\Omega\ww\Omega$
which remains valid restricted to $\cf^1$. 
Thus we calculate that  on $\cf^1$, $0=d\oo a= -\ooo aj\ww\oo j$
which, by the Cartan Lemma, implies that
$\ooo aj=\hhh ajk\oo k$ for some functions $\hhh a jk=\hhh a kj$ well
defined
on $\cf^1$. Although these functions vary in the fiber, the tensor
$$
II_x:= \hhh ajk\oo j\oo k\ot\ee a 
$$
descends to be a well defined element of $S^2T^*_xM\ot N_xM$,
called the {\it second fundamental form} of $M$ at $x$.
$M$ is said to be {\it minimal} if $\ttrace_g(II)=0$, where
$g$ is the induced Riemannian metric on $M$. Such $M$ are locally
of least volume in the sense explained before for surfaces.

\subheading{An exterior differential system
for   submanifolds of Euclidean space}
An eds  for lifts of submanifolds in $\bee{n+s}$ may be defined on
  $  \cf \times \brr{\binom{n+1}2s}$
where $\brr{\binom{n+1}2s }$ has coordinates $h^a_{ij}=h^a_{ji}$,
$1\leq i,j\leq n,\ 3\leq a,b\leq 2+s$   The system is
$$
\align &
\ci = \{ \oo a, \taa a i:= \ooo ai-\hhh aij\oo j\rdiff \tag sub.1\\
&\Omega = \oo 1\ww\hdots\ww \oo n\endalign
$$
This is the analog of the cannonical contact system on $J^1$.
To study submanifolds satisfying a differential equation,
we   restrict to the submanifold of $\cf\times \brr{\binom{n+1}2s}$
where the equation is satisfied.

\subheading{An exterior differential system
for minimal surfaces}
An eds  for lifts of minimal surfaces in $\bee{2+s}$ may be defined on
  $\Sigma\subset \cf \times \brr{3s}$
where $\brr{3s}$ has coordinates $h^a_{ij}=h^a_{ji}$,
$1\leq i,j\leq 2,\ 3\leq a,b\leq 2+s$ and  where
$\Sigma$ is defined by the $s$ equations $\hhh a11+\hhh a22 =0$.

\smallpagebreak
 
\noindent{\bf Remark}. Note that the condition $h_{11}+h_{22}$ is 
significantly simpler
than the equation for minimal surfaces
presented in example 1. This illustrates  the advantage of frames.

\smallpagebreak

Our  system is
the restriction of (sub.1) to $\Sigma$.
Differentiating, we have
$$
\align
d\oo a& \equiv 0\tmod \{ \oo a, \taa a i \} \tag sub.2\\
d\pmatrix
 \taa a 1\\ \taa a 2\endpmatrix &\equiv 
\pmatrix
2\hhh a21\ooo 21 - d\hhh a11 & -(2\hhh a11\ooo 21- d\hhh a 12)\\
-(2\hhh a11\ooo 21- d\hhh a 12)& -(2\hhh a21\ooo 21 - d\hhh a11)
\endpmatrix
\ww\pmatrix\oo 1\\ \oo 2\endpmatrix \tmod \{ \oo a, \taa a i \} \\
&\equiv \pmatrix  \pii a1 & \pii a2 \\ \pii a2 & - \pii a1\endpmatrix
\ww\pmatrix\oo 1\\ \oo 2\endpmatrix \tmod \{ \oo a, \taa a i \} 
\endalign
$$
where the last line defines the one-forms $\pii a1,\pii a2$.

Here there are Cauchy characteristic vector fields.
  The distribution they span is
$
\{\oo i,\ta a,\pii ai \}\upperp
$. (\lq\lq $\upperp$\rq\rq\  refers to the annhilator in the dual space.)

If we reverse the order of the equations in (sub.2) and write
$$
 d\pmatrix
\taa a 2\\ \taa a 1\endpmatrix \equiv 
  \pmatrix  \pii a2 & - \pii a1 \\ \pii a1 & \pii a2\endpmatrix
\ww\pmatrix\oo 1\\ \oo 2\endpmatrix \tmod \{ \oo a, \taa a i \} 
$$
the system
\lq looks like\rq\ the Cauchy-Riemann equations as presented in example 2. 

On $T^{*\bcc{}}\Sigma$ let
$$
\omega := \oo 1 + i\oo 2,\ \ta a :
= \taa a 2 + i \taa a 1 ,\ \pi^a := \pii a 1 + i\pii a 2
$$

Noting that 
$$
\align
d\omega &= i\ooo 21\ww\omega\\
d\theta & = -\pi\ww\omega\\
d\pi & = 2\pi\ww\ooo 21 + 2(h_{11}+ih_{21})K\omega\ww\overline\omega
\endalign
$$
we see the almost complex structure is integrable, and moreover
  the second equation implies $d\theta\ww\theta\neq 0$.
Thus we have a complex contact structure.
From this we recover the result that minimal surfaces in $\bee{2+s}$
are
given by holomorphic data. To recover an actual coordinate
presentation, one would have to work more.  
In summary

\proclaim{Theorem (see e.g. [Ch, p. 19])} Any minimal surface
$M^2\subset\bee{2+s}$ is the
transform of a holomorphic contact curve.\endproclaim

There is no nice description of all local solutions to
the minimal submanifold system for $n>2$. To understand the situation
better,
it is often useful to look  at special classes of solutions,
e.g. solutions with symmetry. (see e.g. [HsL]). 
Taking a broad view
of what consititutes \lq\lq symmetry\rq\rq ,
it has something to do with reducing the complexity of the problem
via the action of some group. 
For example, if $n$ and $s$ are  even,
a natural sub-class of minimal submanifolds of $\bee\ns$ are the
complex submanifolds (we will see in a minute that these are all
minimal). I will  describe 
some special classes of minimal submanifolds shortly, but first
I'll take a short detour to discuss minimizing submanifolds.

\subheading{ Calibrated submanifolds (Harvey-Lawson [HaL])}

Let $(X^{\ns}, g)$ be a Riemannian manifold. A submanifold
$M^n\subset X$ is said to be {\it minimizing} if it is minimal and
such that for all $x\in M$, $\tilde U_x=M$, where $\tilde U_x$ is the open
neighborhood described in example 1.
While being minimal is a local property, being minmizing is a global
matter and in general is difficult to prove. Harvey and Lawson,
building on Wirtinger's and Federer's  observations about the K\"ahler
form, made the following  definition:

Let $(X^{n+s},g)$ be a Riemannian manifold.
A {\it calibration} $\phi \in\Omega^n(X)$ is
  an $n$-form   having the properties:

i. $d\phi =0$ 

 ii. $\phi_x(E)\leq 1$ for all $E\in G(n,T_xX)$,
for all $x\in X$ (where here  $G(n,T_xX)$ denotes the Grassmanian of
unit $n$-planes).  

\proclaim{Theorem [HaL]}Let $(X^{\ns}, g)$ be a Riemannian manifold
with a calibration $\phi\in\Omega^n(X)$.
 Let $i: M^n\ra X^{\ns}$ be such that $i^*(\phi )
=dvol_M$. Then $M$ is minimizing.\endproclaim
\demo{Proof}
Let $U\subset M$ be any open neighborhood and let $V^n\subset X$ be
such that 
$\partial U=\partial V$. Then we have
$$
\text{volume}( U)= \int_U\phi = \int_V\phi \leq \text{volume}(V) 
$$
where the second equality is due to Stokes' theorem.\qed\enddemo

Given a calibration $\phi$, let $face_x(\phi )\subset G(n,T_xX)=
\{ E\in G(n,T_xX) | \phi (E)=1\}$. If $X= \Bbb E^{\ns}$
and $\phi$ has constant coefficients, one can talk about the
face of $\phi$ in $G(n,\ns)$. (Here, and in what follows,
$G(n,\ns)$ denotes the orthogonal Grassmanian.) Note that if
$\gamma (M)\subset face(\phi )$ then
$M$ is minimizing.

In summary:
$$
\matrix 
M \text{ is minimal}&
\leftrightarrow &
 \text{conditions on the}&
&
\leftrightarrow &\text{second order system of pde}\\
& & 
 \text{derivative of }\gamma
& & & \\
M \text{ is minimizing} &\leftarrow &
 \text{conditions on  the  }   &
& \leftrightarrow &\text{first order system of pde}\\
& & 
 \text{image of }\gamma
& & &  
\endmatrix
$$

{\bf Examples}:
In [HaL], Harvey and Lawson  study the following  faces of calibrations:

\noindent 0. The $SU(m+r)$ orbit of a complex $m$-plane, which is
the
complex Grassmannian $G(\bold C^m, \bold C^{m+r})\subset  G_{2m,2(m+r)}$.
(Already studied by Wirtinger and Federer.)

\noindent 1. The $SU(n)$ orbit of a real $n$-plane
in $\bee {2n}$, called the {\it special Lagrangian face}
$\subset G(n,2n)$.

\noindent 2. The $G_2$ orbit of an associative
 $3$-plane  in $\bee 7$, called the
{\it associative  face} $\subset G(3,7)$.  

\noindent 3. The $G_2$ orbit of a coassociative 4-plane in $\bee 7$, called
the
 {\it coassociative face} $\subset G(4,7)$.  

\noindent 4. The $Spin(7)$ orbit of a quaternionic $4$-plane in $\bee 8$,
called the {\it Cayley face} $\subset G(4,8)$.

\smallpagebreak

Harvey and Lawson proved all these faces are
{\it involutive} in the sense of the Cartan-K\"ahler theorem,
that is, that they admit a large class of minimizing submanifolds.

\smallpagebreak

There are not many examples or classes of examples of
minimal submanifolds studied explicitly. To produce some new classes of
examples,
I made the following definition, generalizing the notion of the
faces of calibrations:

\noindent{\bf Definition}. An {\it m-subset} is a submanifold
$\Sigma\subset
G(n,n+s)$ with the property that any $M^n\subset\bee\ns$ with
$\gamma (M)\subset\Sigma$ is necessarily a minimal submanifold.
Given an $m$-subset $\Sigma$, we call any $M$ with $\gamma (M)\subset
\Sigma$
a {\it $\Sigma$-manifold}.

An $m$-subset $\Sigma$ corresponds to a
  first order system of pde that implies the minimal
submanifold system. What follows is a description of an
$m$-subset that is not a face:
 
If we consider $\bee 5$ as the space of $3\times 3$ traceless-symmetric
matrices, then $SO(3)$ acts on $\bee 5$ by conjugation (i.e. given
$x\in \bee 5$ and $g\in SO(3)$, the action is $x\mapsto gxg\inv$).
Let $\rho :SO(3)\ra SO(5)$ denote the representation.
We will use this $\rho (SO(3))$
action as our symmetry. Since $\rho (SO(3))$ acts on
$\bee 5$, it acts on the linear subspaces, in particular on 
the Grassmannian
$G(3,5)$ (here and in what follows I mean the orthogonal Grassmanian).
There are two special orbits
 of dimension two (where a maximal torus  preserves the $3$-plane).
Let $\Sigma\subset G(3,5)$ be one such orbit.
(Notational apology- this $\Sigma$ has nothing to do
with the spaces we define the exterior
differential systems on.) When one sets up the corresponding
eds one finds the following:

\proclaim{Theorem [L1]} $\Sigma$ is an $m$-subset, that is
any $M^3\subset\bee 5$ with $\gamma (M)\subset\Sigma$ is
minimal.

  Let  $\Cal L$ denote
the space of all
$\rho (SO(3))$- weight zero lines in $\bee 5$ (that is, the weight zero
lines through the origin and their translates). Then $\Cal L$ is a complex
contact manifold, and $\Sigma$-manifolds
are exactly the transforms of the complex contact curves in $\Cal L$.
\endproclaim

One may also think of the weight zero lines through the origin as the
lines through rank one matrices.

How did I find the space $\Cal L$?
$\Cal L$ corresponds to $\cf$ quotiented by the Cauchy characteristics.

The picture is as follows:

$$\Cal F$$

\bigpagebreak

$$\ \Cal L \ \ \ \ \ \ \ \ \ \ \ \bee 5
$$

Note that the transform of a point in $\Cal L$ is a line in $\bee 5$.
One can choose coordinates on $\Cal L$ and compatible coordinates
on $\bee 5$ to get an explicit Weierstrass-type presentation formula:

\proclaim{Theorem [L1, 5.1]} Given a holomorphic function of one variable
$h(z)$, setting $w=h(z), \,\, y=\frac {dh}{dz}$ one obtains a minimal
$3$-manifold in
$\bee 5$ given in coordinates by
$$
x^0=
\frac {1-4|z|^2 + |z|^4 }{2\sqrt 3(1+|z|^2)^2}t   
+\frac{2(-2 + 2|z|^2 + |z|^4)}
{\sqrt 3(1+|z|^2)^4}Re(\bar z^2w) 
-\frac{-5 + 2|z|^2 + |z|^4}{2\sqrt 3(1+|z|^2)^3}Re(\bar z y)
 $$
$$x^1 +i x^2 =
\frac{z(1- |z|^2)}{(1+|z|^2)^2}t 
-\frac{2\bar z |z|^2(2  + |z|^2)}{(1+|z|^2)^4}w -
\frac{2z^3}{(1+|z|^2)^4}\bar w -
\frac{1 -2|z|^2 - |z|^4}{2(1+|z|^2)^3}y
+ \frac{z^2}{(1+|z|^2)^3}\bar y
$$
$$x^3 + ix^4 =
\frac{z^2}{(1+|z|^2)^2}t +
\frac{1 + 4|z|^2 + 2|z|^4}{(1+|z|^2)^4}w -
\frac{z^4}{(1+|z|^2)^4}\bar w
-\frac{z(|z|^2 +2)}{2(1+|z|^2)^3}y
+\frac{z^3}{2(1+|z|^2)^3}\bar{y}
 $$
\endproclaim

\noindent{\bf Remark}:
  These are not all the solutions, only solutions having invertible Gauss
maps
are obtained.  For example $3$-planes are not among the above solutions. 
One
could derive a more general formula, but it would involve integrals,
as in the classical formula. (Conversely, one can write a Weierstrass
type formula in the classical case without an integral for a large class
of minimal surfaces in $\bee 3$.)

One could ask about the simplest such solution, the transform of
the function $h(z)\equiv 0$.

\proclaim{Proposition [L2, 5.2]} \  The minimal submanifold given by
$h(z)\equiv 0$ is the cone over the 
real Veronese, that is,
the set of rank one matrices
whose repeated  eigenvalue is positive. 
\endproclaim

Some more examples of $m$-subsets:

\proclaim{Theorem [L1,2.1]} Let $\rho: SO(3)\ra SO(2n+1)$ be 
an irreducible representation and let $\Sigma\subset G(k, 2n+1)$
be any two dimensional $\rho (SO(3))$ orbit. Then
$\Sigma$ is an $m$-subset.\endproclaim

Recalling that $\frak s\frak o (3) = \frak s\frak u (2)$
another generalization is to $SU(m)$ orbits:

\proclaim{Theorem [L2, 3.10]}
Let $V^N$ denote any irreducible real   $SU(n+1)$  module other
than the adjoint.  Then the orbits 
$\Sigma$ of certain codimension $2$-planes
 are m-subsets.  
 
The $\Sigma$-manifolds may be constructed 
in a manner analgous to the situation above.
One constructs a space of special linear spaces,  $\Cal L$, which is 
  is a complex contact manifold. 
The $\Sigma$-manifolds
 are the transforms  complex $n$-folds in $\Cal L$ that
solve the differential system induced from $\Cal F$. 

   The $\Sigma$-manifolds corresponding to the simplest solutions of the
differential system on $\Cal L$ are homogeneous cones,
 in fact the $SU(n+1)$ orbits
of  special $(N-2n-2)$-planes in $V$.  
\endproclaim

One can also construct inhomogeneous examples:

Recall that $G(n,n+2)$ and $G(2,2+s)$ are complex manifolds, in
fact   quadric hypersurfaces $Q^n\subset\Bbb C\Bbb P^{n+1}$,
$Q^s\subset\Bbb C\Bbb P^{s+1}$.

  Bryant (personal communication)   proved  
that 
complex submanifolds $\Sigma\subset G(2,2+s)$ whose tangent
spaces contain no decomposable vectors are m-subsets.
 (While it was known classically
that a surface was minimal iff its Gauss map was holomorphic,
Bryant proved it was sufficient that the Gauss image was a complex
submanifold satisfying a transversality condition.)

\proclaim{Proposition [L2, 3.8]}
Complex submanifolds $\Sigma\subset G(n,n+2)$ whose tangent
spaces contain no decomposable vectors are m-subsets.
\endproclaim 

\proclaim {Corollary [L2, 3.9]} The generic
 complex submanifold of $G(n,n+2)$ of 
complex dimension $[\frac{n}{2}]$
is an m-subset.
\endproclaim

  One might hope to get new $m$-subsets by deforming faces. 
One can rephrase [L2, 3.8] as:

\proclaim{Theorem [L2]} The complex
 Grassmanians  $G(\bcc{m-1},\bcc m)\subset G(2(m-1), 2m)$ admit
deformations
in the category of m-subsets.
\endproclaim

A deeper result is the following:

\proclaim{Theorem [L2, 3.1,3.1${}^*$]} All the faces 0-4 studied
by Harvey and Lawson except $G( \bcc 1,\bcc m)$ and  $G(\bcc{m-1},\bcc m)$
are rigid in the larger category of $m$-subsets.
\endproclaim

Note that the deformation theory was not previously known even 
in the restricted category of faces.

\smallpagebreak

Because it will be useful for our next lecture, let me now explain
how one looks for $m$-subsets.

Recall that $T_EG(n,\ns)\simeq E^*\ot E\upperp$. 
Given a submanifold
$i: M^n\ra\bee\ns$, we have its {\it Gauss map}
$\gamma : M\ra G(n,\ns )$, where $x\mapsto T_xM$.
Consider
the derivative of the
Gauss map, $ d\gamma_x\in
T^*_xM\ot T_{T_xM}G(n,\ns)=
T^*_xM\ot (T^*_xM\ot N_xM)$. 
Now $d\gamma_x\simeq II_x$, and since the second
fundamental form is symmetric, we have
$
d\gamma_x \in (S^2T^*_xM\ot N_xM)$. 

Now if $\Sigma\subset G(n,\ns )$ is any submanifold, then
$T_E\Sigma\subset E^*\ot E\upperp$ is some linear subspace.
If $\phi : M\ra \Sigma$ is any map, then $d\phi_x\in T^*_xM\ot
T_{\phi (x)}\Sigma$, so
 if
$\gamma (M)\subset\Sigma$, then we must have for all $x\in M$ that
$$
d\gamma_x \in (S^2T^*_xM\ot N_xM)\cap ( T^*_xM\ot T_{T_xM}\Sigma). \tag
$\gamma$
$$
where the first condition is because $d\gamma_x$ is a second fundamental
form and the second because the image lies in $\Sigma$.

More generally, given a linear subspace $A\subset V^*\ot W$,
define the {\it prolongation} of $A$, $A\up 1$, by
$$
A\up 1 := (A\ot V^*)\cap (W\ot S^2V^*)
$$
We may think of of $A$ as specifiying a first order
constant coefficient  system of pde
for mappings $f:V\ra W$, namely the system whose solutions
are those mappings $f$ such that $Jac (f)_x\in A$ for all $x\in V$.
In other words, $A$ consists of the admissible first order
terms for the Taylor series of $f$.
With this perspective,
$A\up 1$ has the interpretation of the admissible 
second order terms in the Taylor series. We will see that $A\up 1$
also may be identified with the space of integral elements
of the corresponding eds
at each point. The condition $(\gamma )$ is that
$d\gamma_x\in (T_{T_xM}\Sigma )\up 1$.

\proclaim{The fundamental observation for the study
of $m$-subsets [L2]}
If $\Sigma\subset G(n,\ns )$ is such that for all $E\in\Sigma$,
$$
(T_E\Sigma )\up 1 \subseteq S^2E^*_0\ot E\upperp
$$
where $S^2E^*_0$ denotes the traceless elements,
then $\Sigma$ is an $m$-subset.
\endproclaim

Thus the infinitesimal
version of the problem to find
examples of $m$-subsets is to
find first order constant coefficient systems that imply the Laplace
system.
(In practice, one restricts to {\it involutive} systems, as will be defined
in the next lecture.)

The theorems above were proved
first by studying the linear problem,
that is studying first order constant coefficient
systems of pde that imply Laplace's system
and then by setting up the appropriate exterior
differential system for $m$-subsets whose tangent spaces (linearized 
system of pde) correspond to solutions of the linearized system.

\subheading{Problems for lecture 2}  

1. Define a canonical system for all $n$-dimensional submanifolds
of $\bee \ns$ on $\cf$. Show that the system (sub.1) is
the prolongation of this system. What are the Cauchy Characteristic
vector fields of this system? Give a geometric interpretation
of the quotient space $\cf /\text{leaves}$.

2. Write down an eds for surfaces in $\bee 3$ with Gauss curvature
identically one. What is the dimension of the space of integral elements
at a point?

3.  Write down  an eds for surfaces
of revolution in $\bee 3$ with Gauss curvature
identically one. Find all
integral manifolds of the system.

4. Show that the eds for constant mean curvature one surfaces in hyperbolic
three space is equivalent to the Cauchy-Riemann system. What about
constant mean curvature in $\Bbb E^3, S^3$?

\heading Lecture 3: The Cartan algorithm via the isometric embedding
problem
\endheading

\proclaim{Example 4, The isometric embedding problem}
Let $(M^n,g)$ be a patch of an analytic Riemannian manifold.
What (if any) are the isometric embeddings $M^n\ra \bee\nr$?
\endproclaim

The question is stated a little vaguely- we first want to know
if, given $r$, there exist any embeddings, and then, if so,
\lq\lq how many\rq\rq ?

\def\cfr{\cf_{\bee\nr}}
\def\cfm{\cf_{M}}
\def\et#1{\eta^{#1}}
\def\ett#1#2{\eta^{#1}_{#2}}

To begin, let $\cf_M$ and $\cf_{\bee\nr}$
be the respective frame bundles with $1\leq i,j\leq n$,
$n+1\leq\mu,\nu\leq\nr$. We use $\eta$ to denote forms on $\cf_M $ and
$\omega$ forms on $\cfr$. On 
$\cfm$ we have equations:
$$
\align
dx&= \et i\ee i\\
d\ee i& = \ett ji\ee j\\
d\et j& = -\ett jk\ww\et k \\
d\ett ij&= -\ett il\ww\ett lj + \rr ijkl\et k\ww\et l \endalign
$$
where $\rr ijkl$ are well defined functions on $\cfm$, the coefficients
of the Riemann curvature tensor.

If $\Gamma \subset \cfm\times\cfr$ is the graph of an isometric embedding,
it is natural to demand that  on $\Gamma$
 $\tspan \{\oo i\} = \tspan \{ \et i\}$ and
that $\oo\mu =0$. (Here and in what follows,
we commit a standard abuse of notation,  
what  we really mean that if the inclusion map  is
$i :\Gamma\ra \cfm\times\cfr$, that $i^*(\oo\mu )=0$ etc...).
In addition, it is natural to eliminate some of the Cauchy characteristics
by requiring the bases to line up, that is requiring $\oo i=\et i$
on $\Gamma$.
With these requirements, we obtain, on $\cfm\times\cfr$, the  
Pfaffian system with independence condition
$$
\align 
\ci_0 &= \{ \oo\mu ,\oo j -\et j\}  \\
\Omega &= \et 1\ww\hdots\ww\et n\endalign
$$

From experience, it is also clear what to do next, namely differentiate
(again, we have defined a distribution and now we want to
see how the different planes of the distribution \lq fit together\rq).
We find
$$
\align
d\oo \mu &\equiv -\ooo\mu j\ww\et j\tmod \ci_0 \tag iso.1\\
d(\oo i-\et i )&\equiv - (\ooo ij-\ett ij )\ww\et j\tmod \ci_0\tag
iso.2\endalign
$$
Using the independence condition and the Cartan Lemma, we see that
for (iso.1), (iso.2) to hold, there must be functions
$u^i_{jk},\hhh\mu jk$, defined on
$\Gamma$, symmetric in their lower indices, such that, on $\Gamma$,
$$
\align
(\ooo ij-\ett ij ) &\equiv u^i_{jk}\et k \tmod \ci_0\\
\ooo\mu j&\equiv\hhh\mu jk\et k\tmod \ci_0\endalign
$$

In fact $u^i_{jk}=0$ for all $i,j,k$ (exercise - why?). 
So we have $\tdim A\up 1= r\binom{n+1}2$. At this point we need to perform 
Cartan's test to see if our inital value problem is well posed.
Informally, write out tableau as
$$
A= \pmatrix \ooo\mu j \\ \ooo ij-\ett ij\endpmatrix
$$
an $(r+n)\times n$ matrix of one-forms. Note that $\tdim A = nr+\binom n2$.
To perform Cartan's test (as explained below), 
let $L_j\subset V^*$ be a generic linear subspace of dimension $j$
and let $A_k= A\cap (L_{n-k}\ot W)$.
We need to calculate
$\tdim A_k$. If our bases are sufficiently generic, $\tdim A_{n-k}$ is
just the number of linearly indendent forms in the first $k$ columns.
We get $\tdim A_{n-k} = rk+ (n-1)+(n-2) + \hdots + (n-k)$.
Cartan's test says that the initial value problem is
well posed, or that the system is {\it involutive} if
$\tdim A\up 1= \tdim A +\Sigma_{k=1}^{n-1}\tdim A_k$.
Here
$$
\Sigma_k\tdim A_k = r\binom{n+1}2 + n(n-1)+(n-1)(n-2) + \hdots
+ 2(1) >r\binom {n+1}2=\tdim A\up 1
$$
so the system is not involutive. Thus we must prolong, which amounts
to adding the space of integral elements to our ambient manifold
and the   forms $\ooo\mu j- \hhh\mu jk\et k, \ooo ij-\ett ij$
to our ideal. That is let
$$
\tilde\Sigma = \cfm\times\cfr\times \brr{r+\binom {n+1}2}
$$
and let 
$$
\ci =\{   \oo\mu ,\oo j -\et j,
\ooo ij-\ett ij ,\ooo\mu j-\hhh\mu jk\et k
\} 
$$

The $\hhh\mu jk$ correspond to the
 coefficients of the second fundamental
 form of the submanifold, so what we have done here is 
to  
 add second
derivatives as independent variables.  
 We have now rigged things
such that on integral elements
$$
\align
d\oo \mu &\equiv0 \tmod \ci \\
d(\oo i-\et i )&\equiv 0\tmod \ci  \endalign
$$
This is always true when one prolongs, the old forms automatically
will have zero derivative modulo the new ideal.
Now we   must compute the other derivatives.
$$
d (\ooo ij -\ett ij )\equiv
(\rr ijkl - \Sigma_{\mu}(\hhh\mu ik\hhh\mu jl -\hhh\mu il\hhh\mu jk))
\et k\ww\et l   \tmod \ci
$$
Here we are presented with a similar problem as in example 1, the right
hand side involves independent forms in the independence condition.
There is no way to have 
$d (\ooo ij -\ett ij )=0$ on an integral manifold unless
all the functions
$$
\rr ijkl - \Sigma_{\mu}(\hhh\mu ik\hhh\mu jl -\hhh\mu il\hhh\mu jk)\tag
iso.3
$$
are identically zero. In order to ensure this, we start over yet again,
this time on the submanifold
$\Sigma '\subset \Sigma$ defined by requiring the functions in (iso.3)
to be zero. (iso.3) is an example of
the appearance of terms representing {\it torsion}. There is a danger
that when one restricts to the submanifold of the original manifold
defined by setting the torsion equal to zero that  the
independence condition may no longer hold.
When this happens there are no integral manifolds to
the system (as in example 1, case 2b). Fortunately here this is not a
concern, because we
may view the additional equations as cutting down the 
$\brr{r+\binom{n+1}2}$ factor and not touching $\cfm$ in the
product, which is where the independence condition lies.

Now we continue to take derivatives:
$$
d(\ooo \mu j- \hhh\mu jk\oo k) \equiv (
d\hhh\mu jk -\hhh\nu jk\ooo\mu\nu + \hhh \mu jl\ooo l k+
\hhh\mu kl\ooo lj)\ww\oo k \tmod\ci
$$
To simplify notation, let $\pi^{\mu}_{jk}= 
d\hhh\mu jk -\hhh\nu jk\ooo\mu\nu + \hhh \mu jl\ooo l k+
\hhh\mu kl\ooo lj$. Note that
$\Sigma'$ has a coframing by
$ \et i, \ett ij,\oo i,\oo\mu, \ooo\mu j , \pii\mu{jk},\ooo ij,\ooo\mu\nu$.
The last two groups of forms are dual to infinitesimal rotations
in the tangent and normal bundles, and as before, these infinitesimal
rotations are Cauchy characteristic vector fields. That is, the
Cauchy characteristic distribution is
$$
\{  \et i, \ett ij-\ooo ij , \oo i,\oo\mu, \ooo\mu j, \pii\mu{jk}\}\upperp
.
$$

Now we are in a situation where we have relations among differential
forms and their derivatives at a point. 

In the determined case, if one applies the 
Cartan test again  at a general point of $\Sigma '$,
one finds that the system  is in {\it involution}
(as defined below) and one recovers
the BCJS  theorem stated earlier.

In the overdetermined cases, for a general metric there will be torsion and
if one
restricts to the submanifold $\Sigma''\subset \Sigma '$ where the
torsion is zero, one finds that the independence condition
no longer holds so that there are no integral manifolds.
However, for certain special metrics, there are solutions in
some overdetermined cases, e.g. the sphere $S^n\subset\bee{n+1}$.

Here are some additional results on isometric embeddings in
the overdetermined case.  

\proclaim{Theorem (Cartan) [C1,C2]}
Let $(M^n,g)$ be a (patch of a) flat Riemannian manifold,
then   $M$ admits no
linearly full local isometric embedding to $\bee m$ for $n<m<2n$.

$M$ does admit 
linearly full local isometric embeddings to $\bee {2n}$,
and these depend locally on $\binom n2$ functions of two variables.
\endproclaim 

\proclaim{Theorem (Cartan) [C1,C2]}
Let $(M^n,g)$ be a (patch of a)   Riemannian manifold
with constant sectional curvature $-1$,
then  $M$ admits no
  local isometric embedding to $\bee m$ for $m<2n-1$.

$M$ does admit 
 local isometric embeddings to $\bee{ 2n-1}$,
and these depend locally on $n(n-1)$ functions of one variable,
which is the most possible for any Riemannian $n$-fold admitting
local isometric embeddings to $\bee{ 2n-1}$.
\endproclaim 

In [BBG], Cartan's theorem was generalized to a larger class
of Riemannian metrics, called {\it quasi-hyperbolic metrics}. In
[IvL], Ivey and I generalized their results further
to a class of metrics we call {\it quasi-$\kappa$-curved metrics}:

\noindent{\bf Definition}: Let $(M^n,\tilde g)$ be a Riemannian manifold.
We will say $\tilde g$ is a
{\it quasi-$\kappa$-curved metric} if there exists a
smooth positive definite quadratic form $Q$ on $M$
such that for all $x\in M$
$$
R_x=
-\gamma (Q_x,Q_x) +(\kappa +1)\gamma(\tilde g_x,\tilde g_x)
$$
where $\gamma :S^2T^*\ra S^2(\Lambda^2T^*)$
denotes the algebraic Gauss mapping and
$R_x$ the Riemann curvature tensor.  

 Following [BBG] (and
the real work was done by them and Cartan), we showed:

\def\k{\kappa}

\proclaim{Theorem [IvL, A]}  Let $(M^n,g)$, $n\ge 3$, be a
quasi-$\kappa$-curved
 Riemannian manifold with form $Q$.
Let $X^{2n-1}(\k+1)$ be a space form with constant sectional
curvature $\k+1$. Then there exist local isometric embeddings
$M^n \hookrightarrow X^{2n-1}(\k+1)$, with local solutions depending
on $n^2-n$ functions of one variable, if and only if $\nabla Q$
is a symmetric cubic form on $M$ and
$\nabla Q = L\cdot Q$
for some linear form $L\in \Omega^1(M)$.
\endproclaim

In [BBG], quasi-hyperbolic $n$-folds were characterized as 
space-like hypersurfaces
in $\Bbb L^{n+1}$. We observed that quasi-$\kappa$-curved $n$-folds
arise naturally as spacelike hypersurfaces of a Lorentzian 
sphere having  radius $\frac a{\sqrt{\kappa +1}}$ in $\Bbb L^{n+2}$
and we determined which of these admits optimal isometric
embeddings [IvL, theorem B]. 

Our original motivation was   to look for interesting classes
of minimal isometric embeddings, what we actually 
ended up proving were several
non-existence and rigidity results [IvL, theorems C,D,E].

 \bigpagebreak

\heading Cartan's algorithm for linear Pfaffian systems\endheading

Given $(I,J)$, a linear Pfaffian system on $\Sigma$.

\smallpagebreak

\noindent{\bf Step 1}. Take
a coframing of $\Sigma$ adapted to the filtration
$I\subset J\subset T^*\Sigma$ with forms
$(\ta a, \oo i,\pi^{\epsilon })$, which we assume to be generic
among such adapted coframings. Fix $x\in \Sigma$
and write $V^*=(J/I)_x$, $W^*=I_x$,
$v^i=\ooo ix, w^a=\taa a x$ and $v_i,w_a$
the corresponding dual basis vectors.
\smallpagebreak

\noindent{\bf Step 2}. Calculate $d\ta a$. Since the system is linear
Pfaffian, 
$d\ta a$ is of the following form:
$$
d\ta a\equiv 
\aaaa a\epsilon i\pi^\epsilon\ww\oo i + \tor aij\oo i\ww\oo j\tmod I
$$
The forms $\pii ai$ used previously were 
$\pii ai= \aaaa a\epsilon i\pi^{\epsilon}$.
If one wishes to work on a quotient space, the Cauchy characteristic
distribution is $\{ \ta a, \oo i, \pii ai\}\upperp$.

Define the {\it tableau}
$$
A=A_x:=
\{\aaaa a\epsilon iv^i\ot w_a\subseteq V^*\ot W \mid 1\leq\epsilon\leq r\}
$$
Let $\delta$ denote the natural   the skew symmetrization  map 
 $\delta :W\ot V^{* }\ot V^*\ra W\ot \Lambda^2V^*$ and let
$$
 H^{0,2}(A_x):=W\ot \Lambda^2 V^*/\delta (A_x \ot V^* )
$$
The {\it torsion} of $(I,J)$ at $x$ is
$$
[T_x] := 
[T^a_{ij}w_a\ot v^i\ww v^j]\in      H^{0,2}(A_x).
$$
 Note that in our examples the
torsion was (2.2) in example 1, and (iso.3) in example 4.
\smallpagebreak

\noindent{\bf Step 3}. If $[T]\neq 0$, then start again on
$\Sigma'\subset\Sigma$ defined
by the equations   $[T]=0$.
In practice, since one works infintesimally, one only uses
the equations $[dT]=0$ and checks what relations this forces on
the forms one was using.

\smallpagebreak

\noindent{\bf Step 4}. Assuming $[T]=0$, 
Let $A_{n-j}:= A\cap
(\tspan\{ v^1\hd v^j\}\ot W)$. 
Let $A\up 1:= (A\ot V^*)\cap (W\ot S^2V^*)$, the {\it prolongation} of the
tableau $A$.

\proclaim{Proposition (Cartan) [BCG${}^3$, p120]}
$
\tdim A\up 1 \leq \tdim A +\tdim A_1 + \hdots + \tdim A_{n-1}.
$
\endproclaim

We say $A$ is {\it involutive} if equality holds in the proposition.

In practice it is more convienent to formulate Cartan's inequality  
as follows: Let 
the {\it characters} of $A$,
$s_j$, be the integers defined inductively by
$s_1+\hdots + s_j=\tdim A_{n-j}$.
(In practice this is just the number of independent entries appearing in
the first $j$ columns of $A$,
when $A$ is  written as a subspace of the $\tdim V\times\tdim W$ matrices.)
With this formulation, Cartan's proposition   becomes
$$
\tdim A\up 1 \leq s_1+2s_2 +\hdots + ns_n.
$$

\proclaim{The Cartan-K\"ahler theorem for linear Pfaffian systems}
Let $(I,J)$ be a linear Pfaffian system on $\Sigma$ and let
$x\in \Sigma$ be a point where all 
numerical invariants defined above are locally
constant. If $[T_x]=0$
and $A_x$ is involutive, then there are local solutions to an initial
value problem depending on $s_p$ functions of $p$ variables,
where $s_p$ is the last nonzero character.
\endproclaim

The Cartan-K\"ahler theorem is a generalization of the
Cauchy-Kowalevski theorem, and its proof involves reducing
the problem of finding integral manifolds to a series of 
Cauchy problems.

\smallpagebreak

\noindent{\bf Step 5}. If $A$ is not involutive, one needs to start over on
 a larger space.
Invariantly, one 
prolongs as described in lecture 1. In practice this amounts to enlarging
$\Sigma$ to include
the elements of $A\up 1$ as independant variables, and adding
differential forms to the ideal $\taa ai:= \aaaa a\epsilon i\pi^{\epsilon}
- p^a_{ij}\oo j$, where $p^a_{ij}v^iv^j\ot w_a\in A\up 1$.
In our isometric embedding example, we added $\hhh\mu ij$ as new variables
and the forms $\ooo ij-\ett ij,\ooo\mu j-\hhh\mu jk\oo k$ to the ideal.

\medpagebreak

In summary, we have the following flowchart:

$$\CD
  \matrix \text{Input: Linear Pfaffian}\\ \text{system } (I,J)\text{ on
}\Sigma
\endmatrix
 @<<<   @. \matrix\text{Start over on }\Sigma'\\ \text{rename }\Sigma'
\text{ as }\Sigma\endmatrix\\
   @VVV  @.  @AA{\text{no}}A \\
   \text{Is } [T]\equiv 0?   @>>{\text{no}}> 
\matrix \text{Restrict to }\Sigma'\subset \Sigma
 \\ \text{defined by }[T]=0\endmatrix @>>> 
\text{Is } \Omega\mid_{\Sigma'}\equiv 0 ?
\\
  @VV{\text{yes}}V  @.   @VV{\text{yes}}V\\
   \text{Is tableau involutive?} @>>{\text{yes}}> 
\matrix \text{done: you have} \\ \text{local existence}\endmatrix  
 @. \matrix \text{done: there are no } \\ \text{integral mnflds.}\endmatrix
 \\
@VV{\text{no}}V\\
\matrix \text{prolong to get a new}\\ \text{system on }\tilde\Sigma\subset
G(n,T\Sigma)\\
 \text{Rename }
\tilde\Sigma\text { as }\Sigma \\ \text{ and the }\text{cannonical} \\
\text{system as }(I,J) \endmatrix 
\endCD 
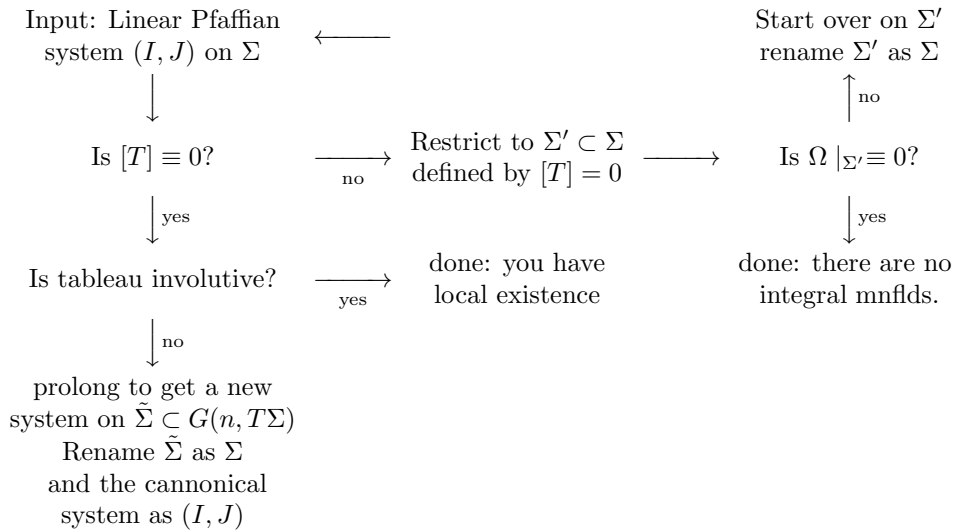
$$

\subheading{Problems for lecture 3 and beyond}

1. Determine the space of integral manifolds for problems 1 and 2
from the problems for lecture 2.

\smallpagebreak

2. Perform Cartan's test for the Cauchy Riemann equations. Interpret
the result in terms of the   holomorphic function $f(z)$.

\smallpagebreak

3. Characterize the surfaces $M^2\subset\bee 3$ such that all points
of $M$ are umbillic points (that is, the two eigenvalues of
the second fundamental form are everywhere equal).

\smallpagebreak

4. Prove that given any nonplanar minimal surface $M^2\subset\bee 3$,
there is  a one-parameter family of  noncongruent deformations that
preserve
the induced Riemannian metric and minimality.

\smallpagebreak

5. Calculate Cartan's test for the case $\Sigma = J^1(\brr n,\brr s)$
with its canonical system and for its prolongation. (Of course the answer
should
be that solutions depend on $s$ functions of $n$ variables.)

\smallpagebreak

6. The Grassmanian $G(n,n+s)$ of oriented orthogonal $n$-planes
in $\brr\ns$ has a natural Riemannian metric on it (such that
the forms $\ooo ai$ form an orthonormal framing). Thus given
$i: M^n\ra\bee\ns$ we can compare the metric induced by $i$ with
that induced by the Gauss map $\gamma$. Suppose that the
two metrics agree pointwise up to a constant, i.e.
$g_i= \lambda g_{\gamma}$, where $\lambda$ is some function.

a. Show that in fact $\lambda$ must be constant.

b. Show that such $M$ are Einstein, that is $Ricci_g= cg$ for
some constant $c$.

c. Classify all such hypersurfaces, i.e.  the case $s=1$. 

d. Classify the case $n=3$, $s=2$. In particular show there
are \lq\lq many\rq\rq\ such. How many?

\smallpagebreak

7. Show that the eds for special Lagrangian submanifolds is involutive.
Integral manifolds depend on how many functions of how many variables?

\smallpagebreak

8. (For exceptional groups enthousiasts) Same question as 6 for
the other faces of calibrations defined by Harvey and Lawson.
 
\subheading{Acknowledgements} I would like to thank the organizers,
especially H. Shin, for their  hospitality and organizing
a great conference. I would also like to thank Daewoo for sponsoring
the conference and their support  of mathematical research.


\Refs

\refstyle{A}
\widestnumber\key{ACGH}


\ref \key BBG \by  E. Berger, R. Bryant, P. Griffiths
\paper  The Gauss equations and
rigidity of
isometric embeddings
\jour  Duke Math. J.
\vol50
\yr 1983
\pages 803-892
\endref

\ref \key BCG${}^3$ \by   Bryant, Chern, Gardner, Goldschmidt, Griffiths
\book  Exterior Differential Systems
\publ Springer Verlag
\publaddr New York
\yr 1991
\endref



\ref \key Ca1 \by  E. Cartan
\paper Sur les vari\'et\'es de courbure constante d'un
espace euclidien ou non euclidien
\jour  Bull. Soc. Math France.
\vol 47
\yr 1919
\pages 125-160
\endref

 \ref \key Ca2 \by  E. Cartan
\paper Sur les vari\'et\'es de courbure constante d'un
espace euclidien ou non euclidien
\jour  Bull. Soc. Math France.
\vol 48
\yr 1920 
\pages 132-208
\endref

  \ref \key Ch \by  S. Chern
\paper  Minimal submanifolds in a Riemannian manifold
\jour  Univ. Kansas Technical report
\vol 19
\yr 1968 
\pages 1-58
\endref

 
 

\ref \key HaL \by R. Harvey and H.B.  Lawson
\paper Calibrated Geometries
\jour Acta Mathematica 
\vol148
\yr 1982
\pages 47-157
\endref


\ref\key HsL\by W.Y. Hsiang and H.B. Lawson
\paper Minimal Submanifolds of Low Cohomogenity
\jour J. Diff. Geom.
\vol  5
\yr1971
\pages 1-38
\endref

\ref \key IlL \by  B. Ilic and J.M. Landsberg
\paper On symmetric degeneracy loci, spaces of symmetric matrices 
of constant rank and dual varieties
\jour  alg-geom/9611025
\endref

 \ref \key IvL \by T. Ivey and J.M. Landsberg
\paper    On minimal isometric embeddings
\jour   Duke
\vol 89 \yr 1997
\pages 1-22
\endref


 
\ref \key L1 \by J.M. Landsberg
\paper   Minimal Submanifolds of $\bold E^{2n+1}$ Arising From
 Degenerate $SO(3)$ Orbits on the Grassmannian
\jour  Trans. A.M.S.
\vol  325 \yr 1991
\pages 101-118
\endref


 
\ref \key L2 \by J.M. Landsberg
\paper   Minimal submanifolds   defined by first order systems
of PDE
\jour  Journal of Differential Geometry
\vol  36 \yr 1992
\pages 369-417
\endref

 
 
\ref \key L3 \by J.M. Landsberg
\paper On second fundamental forms of projective varieties
\jour Inventiones math
\vol 117 \yr 1994
\pages 303--315
\endref


 


  
 \ref \key L4 \by  J.M. Landsberg
\paper     On the local differential
geometry of complete intersections
\jour   S\'eminaire de th\'eorie spectrale et
g\'eom\'etrie, Grenoble 
 \yr 1994-5
\pages 1-12
\endref

 
\ref \key L5 \by J.M. Landsberg
\paper  On degenerate secant and tangential varieties and local
differential  geometry
\jour  Duke Mathematical Journal
\vol 85 \yr 1996
\pages 1-30
\endref
  
 
 \ref \key L6 \by J.M. Landsberg
\paper  Differential-geometric characterizations of complete
intersections
\jour  Journal of Differential Geometry
\vol 44 \yr 1996
\pages 32-73
\endref

\ref \key L7 \by J.M. Landsberg
\paper  Is a linear space contained in a variety? 
- On the number of derivatives needed to tell
\jour  alg-geom 
\endref

\ref \key L8 \by J.M. Landsberg
\paper   On an unusual conjecture of Kontsevich and variants of
Castelnuovo's lemma
\jour  alg-geom/9604023
\endref
 
  


 
\ref \key S \by M. Spivak
\book A comprehensive introduction to differential geometry I-V
\publ Publish or Perish
\publaddr Wilmington
\yr 1979
\endref
 

\endRefs
\enddocument